\documentclass[11pt]{amsart}

\usepackage{amsopn,amsmath,amssymb,eucal,url,enumerate,amscd,amsgen}
\usepackage[pagebackref]{hyperref}

 \newcommand{\nc}{\newcommand}

\nc{\bb}{\mathfrak{b} }
 \nc{\cc}{\mathfrak{c} }  \nc{\dd}{\mathfrak{d} } 
    \nc{\ggo}{\mathfrak{g} }
 \nc{\hh}{\mathfrak{h} }  \nc{\ii}{\mathfrak{i} }
 \nc{\jj}{\mathfrak{j} }  \nc{\kk}{\mathfrak{k} }
\nc{\mm}{\mathfrak{m} }   \nc{\nn}{\mathfrak{n} }
\nc{\pp}{\mathfrak{p} }   
\nc{\rr}{\mathfrak{r} } \nc{\sg}{\mathfrak{s} }
 \nc{\sso}{\mathfrak{so} }  \nc{\spg}{\mathfrak{sp} }
 \nc{\sug}{\mathfrak{su} }  \nc{\ssl}{\mathfrak{sl} }
 \nc{\tog}{\mathfrak{t} }  \nc{\uu}{\mathfrak{u} }
 \nc{\vv}{\mathfrak{v} } \nc{\ww}{\mathfrak{w} }
 \nc{\zz}{\mathfrak{z} }

\nc{\CC}{{\mathbb C}}
 \nc{\DD}{{\mathbb D}}
\nc{\FF}{{\mathbb F}}
\nc{\GG}{{\mathbb G}}  
\nc{\HH}{{\mathbb H}}
\nc{\II}{{\mathbb I}}
\nc{\JJ}{{\mathbb J}}
\nc{\KK}{{\mathbb K}}
\nc{\NN}{{\mathbb N}}

\nc{\RR}{{\mathbb R}}  
 \nc{\ZZ}{{\mathbb Z}}

\nc{\ggob}{\overline{\mathfrak{g}}} 
 
\nc{\glg}{\mathfrak{gl} }
  
\nc{\pca}{\mathcal{P}} \nc{\nca}{\mathcal{N}}
 
 \nc{\vp}{\varphi} \nc{\ddt}{\frac{{\rm d}}{{\rm d}t}}
 \nc{\la}{\langle} \nc{\ra}{\rangle}
 
 \nc{\SO}{{\sf SO}} \nc{\Spe}{{\sf Sp}} \nc{\Sl}{{\sf Sl}}
 \nc{\SU}{{\sf SU}} \nc{\Or}{{\sf O}} \nc{\U}{{\sf U}}
 \nc{\Gl}{{\sf Gl}} \nc{\Se}{{\sf S}} \nc{\Cl}{{\sf Cl}}
 \nc{\Spin}{{\sf Spin}} \nc{\Pin}{{\sf Pin}}
 \nc{\Ta}{{\sf T}}
 
 \nc{\ad}{\operatorname{ad}} \nc{\Ad}{\operatorname{Ad}}
 \nc{\coad}{\operatorname{coad}}
 \nc{\rank}{\operatorname{rank}} \nc{\Irr}{\operatorname{Irr}}
\nc{\Ima}{\operatorname{Im}}
 \nc{\End}{\operatorname{End}} \nc{\Aut}{\operatorname{Aut}}
 \nc{\Inn}{\operatorname{Inn}} \nc{\Der}{\operatorname{Der}}
 \nc{\Ker}{\operatorname{Ker}} \nc{\Iso}{\operatorname{I}}
 \nc{\Le}{\operatorname{L}} \nc{\tr}{\operatorname{tr}}
 \nc{\dif}{\operatorname{d}} \nc{\sen}{\operatorname{sen}}
 \nc{\modu}{\operatorname{mod}} \nc{\Ric}{\operatorname{R}}
 \nc{\Sym}{\operatorname{Sym}} \nc{\sca}{\operatorname{sc}}
 \nc{\scalar}{{\sf s}} \nc{\grad}{\operatorname{grad}}
 \nc{\ricci}{\operatorname{r}} \nc{\riccin}{\operatorname{Ric}}
 \nc{\Lie}{\operatorname{L}} \nc{\ct}{\operatorname{T}}

 \theoremstyle{plain}
 \newtheorem{thm}{Theorem}[section]
 \newtheorem{prop}[thm]{Proposition}
 \newtheorem{cor}[thm]{Corollary}
 \newtheorem{lem}[thm]{Lemma}
 
 \theoremstyle{definition}
 \newtheorem{defn}[thm]{Definition}
 
 \theoremstyle{remark}
 \newtheorem*{rem}{Remark}
 
 \newtheorem{exa}[thm]{Example}
 \newtheorem{exams}[thm]{Examples}

\begin{document}


\title[Extending invariant complex structures]
{Extending invariant complex structures}

\author[R. Campoamor Stursberg]{Rutwig Campoamor Stursberg}
\address{R. Campoamor Stursberg: I.M.I and Depto. Geometría y Topología, Universidad Complutense de Madrid, Spain
}
\email{rutwig@ucm.es}

\author[I. E. Cardoso]{Isolda E. Cardoso}
\address{I. E. Cardoso, ECEN-FCEIA, Universidad Nacional de Rosario \\
Pellegrini 250, 2000 Rosario, Santa Fe, Argentina}
\email{isolda@fceia.unr.edu.ar}

\author[G. P. Ovando]{Gabriela P. Ovando}
\address{G. P. Ovando: CONICET and ECEN-FCEIA, Universidad Nacional de Rosario \\
Pellegrini 250, 2000 Rosario, Santa Fe, Argentina}
\email{gabriela@fceia.unr.edu.ar}


\begin{abstract}   We study the problem of extending a complex structure to a given Lie algebra $\ggo$, which is firstly defined on an ideal $\hh\subset \ggo$. We consider the next situations:   $\hh$ is either complex or it is totally real. The next question is to equip $\ggo$ with  an additional structure, such as a (non)-definite metric or a symplectic structure and to ask either $\hh$ is non-degenerate, isotropic, etc. with respect to this structure, by imposing a compatibility assumption. We show that this implies certain constraints on the algebraic structure  of $\ggo$. Constructive examples illustrating this situation are shown, in particular computations in dimension six are given.
\end{abstract}

\thanks{{\it (2000) Mathematics Subject Classification}: 53C15, 53C55, 53D05, 22E25, 17B56 }


\maketitle

\section{Introduction}

An important source of complex manifolds is provided by homogeneous manifolds $M= G/H$ with trivial  isotropy, that is $H=\{0\}$ so that $M$ is itself a Lie group and the geometric structure is invariant under left-translations. Thus the geometric structure is determined at the Lie algebra level. This setting enables the construction and study of many examples and applications, which in the history (starting by the Erlangen problem by Klein) gave answers to several interesting problems, such as existence of complex, symplectic, pseudo-K\"ahler, K\"ahler but non-symplectic structures, as for instance the Kodaira-Thurston manifold).

 In dimension four a general classification  of  Lie groups provided with invariant complex structures is known  \cite{Ov1,Sn} while  the homogeneous case was more recently completed in \cite{CF} but no general result is known in higher dimensions. This is an open topic of active research now, as for instance the six-dimensional situation starting by  the existence problem of such structures (see \cite{ABD,CO,COUV,EF,FOU} for advances in this direction). Indeed the existence   and the classification problems  become more complicated in higher dimensions. Alternative solutions to this is the consideration of other elements such as an additional geometrical structure related to the (almost) complex structure, giving rise to Hermitian or anti-Hermitian metrics, (pseudo)-K\"ahler or complex symplectic structures, tamed complex structures, etc.

Our approach here is the extension of complex structures on a given Lie algebra $\ggo$. Our motivation and starting point is the following observation: with the exception  of the Lie algebras with Heisenberg commutator (see for instance \cite{Ov2}), {\it  in dimension four, most of the  Lie algebras endowed with a complex structure admit an ideal which is either complex or totally real.} 

\smallskip

This suggests the study of the relationship between the algebraic structure of $\ggo$ and the existence problem of complex structures on $\ggo$ in the following frame: determine a complex structure on $\ggo$ in such way that a fixed ideal $\hh$ is  complex or totally real. 

The  underlying algebraic situation  is the so called {\em extension problem}: extend the structure of the Lie algebra $\hh$ to a Lie algebra $\ggo$ in such way that  $\hh \subseteq \ggo$  is an ideal of $\ggo$, 
 which gives rise to the following exact sequence of Lie algebras:
$$0 \to \hh \to \ggo \to \ggo/\hh \to 0.$$
The extension problem  goes back to Chevalley and Eilenberg \cite{CE} and it is a current research topic in more general situations (see for instance \cite{AM} and references therein).

 As we shall see the existence of a complex or a totally real ideal on $\ggo$  imposes extra conditions on the algebraic structure of $\ggo$. Note that this viewpoint which makes  use of algebraic tools was useful in several works (see \cite{AS,ABD,BF,EFV,FKV,GaV,GR} for instance). 

The next step is to add a non-degenerate symmetric or skew-symmetric bilinear map, that is, a metric or a symplectic structure satisfying some compatibility condition with respect to the complex structure, obtaining Hermitian or anti-Hermitian structures (also known as Norden metrics) or pseudo-K\"ahler or complex symplectic structures. In this case one asks the ideal to be isotropic, non-degenerate, etc. conditions which also impose restrictions on the algebraic structure of $\ggo$. 

Our main results include:
\begin{itemize}
\item the algebraic conditions to extend a complex structure defined on a ideal of a given Lie algebra
 (in terms of representations of Lie algebras and cohomology);
\item the extension of almost Hermitian or almost anti-Hermitian structures starting with an almost Hermitian or almost anti-Hermitian ideal, also symplectic correspondence due to
$$\begin{array}{lll}
\left\{
{\begin{array}{l}
\mbox{Hermitian structures}\\ \mbox{ such that }\nabla J=0 
\end{array}}\right\}
  &   \longleftrightarrow  & 
	\left\{ 
	\begin{array}{l}
\mbox{ pseudo-K\"ahler structures for } J \end{array}\right\}\\
 & &  \\ \\
\left\{
{\begin{array}{l}
\mbox{anti-Hermitian structures}\\ \mbox{ such that } \nabla J=0
\end{array}} \right\}  & \longleftrightarrow  & \left\{\mbox{ complex-symplectic structures  for }J\right\}
\\
 & &  \\
\end{array}
$$
where $\nabla$ denotes the Levi-Civita connection. 

\item Applications in the form of several examples in different situations. In particular computations in dimension six are given providing examples of complex structures on $\ggo$ such that the ideal $\hh$ is totally real.  In this framework this was possible since the integrability condition for the almost complex structure was reduced to a linear equation which was possible to be solved. 
\end{itemize} 


Althought all questions are presented in an algebraic  setting,  several examples and applications are shown along the paper. In this sense  geometrical implications such as curvature, special connections, etc.  have been considered by different authors in particular and separate  cases (see \cite{ABD2,ABDF,AR,CFG}). In the symplectic case, the study of symplectic Lie groups, with special interest on Lagrangian extensions,  was done in \cite{BC}. Our proposal here could be go on in this direction.  

\section{Preliminaries}

In this section we recall basic  facts related to the cohomology of Lie algebras and we introduce  complex structures on Lie algebras in this setting.

\subsection{Cohomology of Lie algebras}

\vskip 3pt

Let $(\kk, [\cdot , \cdot]_{\kk})$ and $(\hh, [\cdot , \cdot]_{\hh})$ be Lie algebras and let $\pi$ denote a representation from $\kk$ into $\hh$ by derivations. 

  Let $C^1(\kk, \pi)$ be the space of linear morphisms of $\kk$ into $\hh$ and for $s >1$ let $C^s(\kk, \pi)$ be the space of $s$-alternanting maps of $\kk \times \hdots \times \kk$ (s factors) into $\hh$.  
  
  The {\em coboundary operator} is a linear operator  $d: C^i(\kk,\pi) \to C^{i+1}(\kk,\pi)$ for all $i\geq 1$. If $\theta$ is an element of $C^1(\kk, \pi)$ one has
\begin{equation}\label{cob1}
d\theta(x,y)=\pi(x) \theta(y)-\pi(y) \theta(x) - \theta([x,y]_{\kk})\qquad \mbox{ for all }x,y \in \kk.
\end{equation}
The map $\theta \to d\theta$ is linear and the {\em 1-cocycles} are the elements  of the kernel of $d$, denoted by  $Z^1(\kk, \pi)$. Let $h \in \hh$ and take $\theta_h\in C^1(\kk,\pi)$ 
\begin{equation}
\theta_h(x)=\pi(x) h \qquad \qquad \mbox{ for } x \in \kk.
\end{equation}
It follows from a trivial calculation that $d\theta_h=0$ for all $h\in \hh$, hence $\theta_h$ is a 1-cocycle which is called a 1-{\em coboundary}; denoting  $B^1(\kk, \pi)$ the set of 1-coboundaries, put 
$$H^1(\kk, \pi)=Z^1(\kk, \pi)/B^1(\kk, \pi).$$
which denotes the first cohomology group of $(\kk, \pi)$.

For $\alpha\in C^2(\kk, \pi)$ let 
\begin{equation}\label{cob2}d\alpha(x,y,z)=  \sum_c \alpha([x,y],z) -\sum_c \pi(x) \alpha(y,z) \qquad \qquad  x,y, z \in \kk
\end{equation}
where $\sum_c$ denotes summation over the set of cyclic permutations of $x,y,z\in \kk$. It is easy to verify that $d\alpha \in C^3(\kk, \pi)$, the map $d:C^2(\kk, \pi) \to C^3(\kk, \pi)$ is linear so that  the set of 2-cocycles is
$$Z^2(\kk, \pi)=\{ \alpha\in C^2(\kk, \pi)\,:\,d\alpha=0\}$$
which is the kernel of $d$. On the other hand if $\theta \in C^1(\kk,\pi)$ then $d^2\theta =0$.

Thus $d$ maps $C^1(\kk,\pi)$ onto a subspace of $C^2(\kk,\pi)$, denoted by $B^2(\kk,\pi)$ and called the 2-coboundaries. Let 
$$H^2(\kk,\pi)=Z^2(\kk, \pi)/B^2(\kk, \pi),$$
be the second cohomology group.

\begin{rem} One extends the operator $d:C^s(\kk,\pi) \to C^{s+1}(\kk,\pi)$  and one defines  the s-cohomology group as
$$H^s(\kk,\pi)=Z^s(\kk,\pi)/B^s(\kk,\pi),$$
 where $Z^s(\kk,\pi)$ denotes the kernel of $d:C^s(\kk,\pi) \to C^{s+1}(\kk,\pi)$  and $B^s(\kk,\pi)=d(C^{s-1}(\kk,\pi))$ denotes the image of $d:C^{s-1}(\kk,\pi) \to C^s(\kk,\pi)$, whose elements are called $s$-coboundaries. 
 \end{rem}

 \begin{exa} Let $\ggo$ be a fixed Lie algebra, and let $\ad$ denote the adjoint representation of $\ggo$, where $\ad(x) y=[x,y]$ for all $x,y \in \ggo$. The Jacobi identity says that the adjoint action acts by derivations of $\ggo$. In the setting above, take $\hh=\kk=\ggo$, so that  $C^1(\ggo, \ad)$ denotes the set of linear morphisms $t:\ggo \to \ggo$ and for $s\geq 2$, let $C^s(\ggo, \ad)=\{ \alpha : \alpha \mbox{ is  s-linear and alternanting on } \ggo \}$.  For $s\geq 1$ the coboundary operator $d:C^s(\ggo, \ad) \to C^{s+1}(\ggo, \ad)$ (see (\ref{cob1}) for $s=1$ and (\ref{cob2}) for $s=2$) induces  the  Chevalley cohomology.
\end{exa}
 
 \vskip 3pt
 
 Let $(\hh, [\cdot, \cdot]_{\hh})$ and $(\kk, [\cdot, \cdot]_{\kk})$ denote
  real Lie algebras, and let $\pi: \kk \to \End(\hh)$ be a representation. 
  Let $\ggo$ be  the vector space direct sum of $\hh$ and $\kk$. 
   We would like to define a Lie algebra structure on  $\ggo=\kk \oplus \hh$ 
   for which  $\hh$ is an ideal. 
 
 Define a skew symmetric bilinear map $[\cdot ,\cdot ]:\ggo \times \ggo \to \ggo$ by
\begin{equation}\label{corchete}
\begin{array}{rcll}
[x,y] & = & [x, y]_{\hh} & x, y\in \hh, \\
{[x,y]} & = & \pi(x) y & x\in \kk, y \in\hh, \\
{[x,y]} & = & [x,y]_{\kk} + \alpha(x,y) & x,y \in \kk,
\end{array}
\end{equation}
where $\alpha: \kk \times \kk \to \hh$ is bilinear  skew-symmetric.

This is a Lie bracket on $\ggo$ if and only if the Jacobi identity holds:
\begin{equation}\label{jacobi}
[[x,y],z]+[[y,z],x] +[[z,x],y]=0\qquad \qquad \mbox{ for all } x,y, z \in \ggo.
\end{equation}
Thus 
\begin{itemize}
\item for $x,y,z \in \hh$, (\ref{jacobi})  follows from the Jacobi identity for 
$[\cdot, \cdot]_{\hh}$; 

\item for $x,y \in \hh$, $z\in \kk$,  (\ref{jacobi}) is satisfied if and only if
$\pi(z)$ is a derivation for every $z\in \kk$:
$$\pi(z)[x,y]_{\hh}= [\pi(z) x, y]_{\hh} + [x, \pi(z) y]_{\hh};$$

\item for $x,y\in \kk$, $z\in \hh$, since $\pi$ is a representation, the Jacobi
identity
reduces to 
$$[\alpha(x,y), z]_{\hh}=0,$$
whence (\ref{jacobi}) is satisfied in this case if and only if  $\alpha$ takes values in the center of $\hh$, denoted by $\zz(\hh)$.

\item for $x,y,z\in \kk$ we have
$$[[x,y],z]= [[x,y]_{\kk} + \alpha(x,y),z]=[[x,y]_{\kk},z]_{\kk} + \alpha([x,y]_{\kk},z) -\pi(z) \alpha(x,y)$$
therefore (\ref{jacobi}) holds if and only if
$$\begin{array}{rcl}
0 & = & [[x,y]_{\kk},z]_{\kk} + [[y,z]_{\kk},x]_{\kk} + [[z,x]_{\kk},y]_{\kk} \\
0 & = & \alpha([x,y]_{\kk},z) - \pi(z) \alpha(x,y)  -\pi(x) \alpha(y,z) + \\
& & + \alpha([y,z]_{\kk},x)  + \alpha([z,x]_{\kk},y) -\pi(y) \alpha(z,x),
 \end{array}
$$
 since  $[\cdot , \cdot]_{\kk}$ is a Lie bracket on $\kk$, the first 
 equality is true. For the second one we ask   $\alpha$ to be  
 a 2-cocycle from $(\kk, \pi)$.
\end{itemize}

The above paragraph proves the following proposition.

\begin{prop}\label{lie} Let $(\hh, [\cdot , \cdot ]_{\hh})$ and 
$(\kk, [\cdot , \cdot ]_{\kk})$ be Lie algebras. Let $\pi$ denote a 
representation from $\kk$ into $\hh$ acting by derivations and let $\alpha \in Z^2(\kk, \pi)$. For $\ggo=\kk \oplus \hh$ direct sum as vector spaces, the bracket $[\cdot , \cdot]: \ggo \times \ggo \to \ggo$  as in (\ref{corchete}) satisfies the Jacobi identity if and only if 
\begin{itemize}
\item the image of $\alpha$ is in the center of $\hh$: $\Ima \alpha \subset  \zz(\hh)$, and
\item  $\alpha \in Z^2(\kk, \pi)$:
$$0 = \alpha([x,y]_{\kk},z) + \alpha([y,z]_{\kk},x)  + \alpha([z,x]_{\kk},y) - \pi(z) \alpha(x,y)  -\pi(x) \alpha(y,z)  -\pi(y) \alpha(z,x).$$
\end{itemize}
\end{prop}

We call the resulting Lie algebra $\ggo=\hh \oplus \kk$ as the {\em extended semidirect product} of $\hh$ and $\kk$ via $(\pi, \alpha)$. Thus one gets the short exact sequence:
$$ 0 \longrightarrow \hh \longrightarrow \ggo \longrightarrow \ggo/\hh \longrightarrow 0.$$

\begin{exams} \label{exas} In the setting above
\begin{itemize}
\item for $\alpha=0$  one gets the semidirect 
product of $\hh$ and $\kk$ via $\pi$. We shall denote the semidirect product as $\ggo=\kk \ltimes_{\pi} \hh$.

In particular the tangent Lie algebra is $\Ta \kk=\kk \ltimes_{\ad} \hh$ where $\hh$ is the underlying vector space to $\kk$ and $\pi=\ad$ the adjoint representation. The cotangent Lie algebra $\Ta^* \kk=\kk\ltimes_{\ad^*} \hh$ where $\hh$ is the underlying vector space to $\kk$ and $\pi=\ad^*$ the coadjoint representation.
\item $\hh$ abelian and $\pi=0$, one gets a central extension of $\kk$;
\item $\hh$ and $\kk$ abelian and $\pi=0$, then $\ggo$ is a 2-step nilpotent Lie algebra.
\end{itemize}
\end{exams}  


\

\begin{prop} \label{p2} Let $(\ggo, [\cdot , \cdot])$ be a Lie algebra which 
decomposes as a direct sum of the vector spaces $\kk$ and $\hh$ where  $\hh$ 
is an ideal in $\ggo$. Then $\kk$ can be endowed with a Lie bracket 
$[\cdot , \cdot]_{\kk}$ in such way that $\ggo=\kk \oplus \hh$ is isomorphic 
to an extended
semidirect product as  constructed above for a suitable pair $(\pi, \alpha)$.
\end{prop} 
\begin{proof} Let $[\cdot, \cdot]_{\hh}:=[\cdot , \cdot]_{\hh \times \hh}$ be the restriction of the Lie bracket of $\ggo$ to $\hh$.  The Jacobi identity in $\ggo$ implies that  $[\cdot, \cdot]_{\hh}$ is a Lie bracket.

Let $\ggo=\kk \oplus \hh$ be a direct sum as vector spaces and let $p_1:\ggo \to \kk$ and $p_2:\ggo \to \hh$ denote the linear projections with respect to the splitting $\ggo=\kk \oplus \hh$.

Since $\hh$ is an ideal in $\ggo$, the quotient space  $\ggo/\hh$ has a Lie algebra structure with Lie bracket $[\cdot, \cdot]'$ for which the projection $p:\ggo \to \ggo/\hh$ is a homomorphism of Lie algebras: $[p x, p y]'= p [x,y ]$ for all $x,y \in \ggo$. Furthermore the restriction of $p$ to $\kk$, $p:\kk \to \ggo/\hh$ is a linear isomorphism and for $x,y \in \kk$ one has $p[x,y]=p p_1 [x,y]$ (identifying $\kk$ with the set $\{(x,0)\in \ggo\,:\, x \in \kk\}$.

 Let $[\cdot, \cdot]_{\kk}$ be the skew symmetric bilinear form on $\kk$ given by $[x,y]_{\kk}=p^{-1}[px, py]'$; it satisfies the Jacobi identity since it is the translation of the Lie bracket $[\cdot, \cdot]'$ by way of $p$. Thus $(\kk, [\cdot, \cdot]_{\kk})$ and $(\ggo/\hh, [\cdot, \cdot]')$ are isomorphic as Lie algebras. Let $\alpha:\kk \times \kk \to \hh$ defined by $\alpha(x,y):=p_2([x,y])$, this is clearly bilinear and skew symmetric. Moreover (via identifications) the relation
 $$[x,y]=[x,y]_{\kk} + \alpha(x,y)$$
holds. Finally for $x\in \kk$, $y\in \hh$ define $\pi(x) y:=[x,y]$.
 
 For $x,y\in \kk$, $z\in \hh$, the Jacobi identity in $\ggo$ implies that $\pi:\kk \to \End(\hh)$ is a homomorphism of Lie algebras and $\Ima \alpha \subset \zz(\hh)$. 
 
 For $x,y\in \hh$, $z\in \kk$ one has $[z, [x,y]]=\pi(z) [x,y]_{\hh}$, therefore the Jacobi identity for $[\cdot, \cdot]$ says that $\pi$ acts by derivations.
 
 For $x,y, z \in \kk$  we have $[[x,y],z]=[[x,y]_{\kk},z]_{\kk}+ \alpha([x,y]_{\kk},z)-\pi(z)\alpha(x,y)$. Hence the Jacobi identity on $\ggo$ implies that $\alpha$ is a 2-cocycle $\alpha \in Z^2(\kk, \pi)$. Finally let $\tilde{\ggo}$ denote the Lie algebra constructed from $(\kk, [\cdot,\cdot]_{\kk})$ and $(\hh,[\cdot,\cdot]_{\hh})$ attached to $(\pi, \alpha)$ then $i_1 + i_2:\tilde{\ggo} \to \ggo$, given by $(i_1 + i_2)(x,y)=x+y \in \ggo$ is an isomorphism of Lie algebras, in fact:
 
 $$\begin{array}{rcl}
 (i_1 + i_2)[(x_1, y_1), (x_2,y_2)]& = & (i_1 + i_2)([x_1, x_2]_{\kk}, [y_1, y_2]_{\hh}+\alpha(x_1, x_2)+\\ && + \pi(x_1) y_2 - \pi(x_2) y_1) \\
 & = & [x_1, x_2]_{\kk}+ [y_1, y_2]_{\hh}+\alpha(x_1, x_2)+ \pi(x_1) y_2 - \pi(x_2) y_1 \\
 & = & [x_1+ y_1, x_2 + y_2].
 \end{array}
 $$ 
 \end{proof}
 
 \begin{rem} \label{prods} It turns out that the structure above describes an {\em almost product structure} on $\ggo$: a linear endomorphism $E:\ggo \to \ggo$ satisfying $E^2=1$ (and $E\neq 1$). Indeed $E$ can be described in terms of its eigenspaces: $\ggo=\ggo_+ \oplus \ggo_-$ where $E_|{_{\ggo_+}}=1$ and $E_|{_{\ggo_-}}=-1$.

When $\dim \ggo_+=\dim \ggo_-$ the almost product structure $E$ is called an {\em almost paracomplex structure}. 

The almost product structure $E$ is {\em integrable} if
$$[Ex,Ey]=-[x,y]+E[Ex,y]+E[x,Ey]\qquad \mbox{ for all }x,y\in\ggo$$
or equivalently the subspaces $ \ggo_+,\ggo_-$ are subalgebras. In this situation $E$ is called a {\em product structure} or a {\em paracomplex structure}. See for instance \cite{ABDO,CFG} and references therein. 
\end{rem} 

 \section{Complex structures on Lie algebras and ideals}
 
 In this section we study the extension problem of complex structures attached to ideals. We consider two situations: the ideal is either invariant by the complex structure or it is totally real.
 
 \medspace
  
An {\it almost complex}  structure on the Lie algebra $\ggo$ is an
element $J\in V^1(\ggo,\ad)$ satisfying $J^2=-1$ (where $1$ is
the identity map). The Nijenhuis tensor for $J$ is defined as
\begin{equation}\label{nijenhuis}
N_J(x,y)=[Jx,Jy]- J d J (x,y) \qquad x,y \in \ggo
\end{equation}
Any almost complex structure $J$ is called {\em integrable} if $N_J\equiv 0$, that is 
\begin{equation}\label{int}
d(J)(x,y) = J^{-1}[Jx,Jy]\qquad \mbox{ for all } x,y \in \ggo,
\end{equation}
or explicitly $[Jx,Jy]=[x,y]+J[Jx,y]+J[x,Jy]$. 


The Nijenhuis tensor verifies, for all $x,y \in \ggo$ the following identities: 

 $$N_J(y,x)=-N_J(x,y)=-N_J(Jx,Jy)\qquad N_J(Jx,y)=N_J(x,Jy)=-JN_J(x,y).$$

Hence if $\ggo$ decomposes into a direct sum of vector subspaces $\ggo= \uu \oplus J\uu$, then  $N_J\equiv 0$   if and only if $N_J(x,y)=0$ for all $x,y \in \uu$. As usual, integrable almost complex structures  are called {\it complex structures}.

Almost complex structures  $J :\ggo \to \ggo$ satisfying one of
the following conditions for any $x, y\in \ggo$:
$$\begin{array}{lrcll}
\mbox{c1)} & \,\,J[x,y] & = & [x,Jy] & 
\\
\mbox{c2)} & \,\,[Jx,Jy] & = & [x,y] & 
\end{array}
$$  are always integrable. 
Complex structures of type c1) determine a structure of  complex
Lie algebra on $\ggo$, they are sometimes called {\em bi-invariant}.
Structures of type c2) are called {\em abelian}.

One has the following  equivalence relation between Lie algebras with complex
structures.  Lie algebras with complex structures $(\ggo_1,J_1)$
and $(\ggo_2,J_2)$ are called {\em equivalent} if there exists an
isomorphism of Lie algebras $\sigma:\ggo_1 \to \ggo_2$ such that
$J_2 \circ \sigma = \sigma \circ J_1$. In particular when
$\ggo_1=\ggo_2$ a classification of complex structures can be
done. See \cite{Sn, Ov1} for the classification in dimension four.

If $J'=\sigma J \sigma^{-1}$, by using that
$\sigma$ is an automorphism one gets that $J'$ is abelian (bi-invariant) if
 $J$ is of this type. 

Let $\vv\subseteq \ggo$ be a subspace on a Lie algebra $\ggo$ equipped
with a complex structure $J$, recall that  $\vv$ is called
$$ \begin{array}{ll}
\mbox{ complex if}\quad & J\vv \subseteq \vv, \qquad \\
\mbox{ totally real if } \quad &
\vv \cap J\vv=\{0\}.
\end{array}
$$

 Now we are interested in studying complex structures  on extended semidirect 
 products  $\ggo=\kk \oplus \hh$ attached to $(\pi, \alpha)$, specifically
 when  the ideal $\hh$ is either complex or totally real.

\

{\bf A complex ideal.}  Let $\ggo$ be a Lie algebra such that $\hh \subset \ggo$ is an ideal. Let $\kk\subset \ggo$ be a complementary subspace of $\ggo$ attached with the pair $(\pi,\alpha)$ and endowed with the algebraic structure given in Proposition \ref{p2}. 

Assume there is a complex structure $J$ on $\ggo$ such that $\hh$ is $J$-invariant. In terms of the direct sum as vector spaces $\ggo=\kk \oplus \hh$, we notice that  the subspace $\kk$ is not necessarily $J$-invariant.  Thus for $x\in \kk$ one has:
$$
J(x)= j(x) + \beta(x)\quad \mbox{ where } j: \kk \to \kk \mbox{ and } \beta: \kk \to \hh \mbox{ is linear}.$$
Since $J^2 = -1$ one gets
$$J^2(x) = -x = j^2(x)+\beta(j(x)) + J\beta(x)$$
and this implies 
$$j: \kk \to \kk \mbox{ is an almost complex structure  and  }$$
\begin{equation}\label{jbeta}
 J\beta (x)=-\beta(j(x)) \quad \forall x\in \kk.
 \end{equation}

 Now the integrability condition of $J$ says:
\begin{itemize}
\item For $x,y\in \hh$
$$[Jx, Jy]_{\hh}  =  [x,y]_{\hh} + J[Jx,y]_{\hh}+J[x,Jy]_{\hh}$$
which is the integrability condition for the restriction  $J_{|_{\hh}}:\hh \to \hh$.
 
\item For $x\in \kk, y\in \hh$, on the one hand
$$[Jx, Jy]= [j(x) + \beta(x), Jy]=\pi(j(x)) Jy + [\beta(x), Jy]_{\hh}$$
and on the other hand
$$\begin{array}{rcl}
[x,y]+ J[Jx,y]+J[x,Jy] & = & \pi(x) y+ J[j(x)+\beta(x), y]+ J[x, Jy]\\
& = & \pi(x) y + J\pi(j(x))y + J[\beta(x), y]_{\hh}+ J \pi(x)Jy ,
\end{array}
$$
hence
\begin{equation}\label{jkh}
\pi(j(x)) Jy + [\beta(x), Jy]_{\hh} = \pi(x) y + J\pi(j(x))y + J\pi(x)Jy  + J[\beta(x), y]_{\hh}.
\end{equation}

\item For $x,y\in \kk$: on the one hand
$$
\begin{array}{rcl}
[Jx,Jy] & = & [j(x)+\beta(x),jy+\beta(y)] \\
 & = & [j(x), j(y)]_{\kk}+ \alpha(j(x),j(y))+\pi(j(x))\beta (y)- \pi(j(y)) \beta(x) + [\beta(x), \beta(y)]_{\hh}
\end{array}
$$ while on the other side
$$
\begin{array}{rcl}
[Jx,Jy]
& = & [x,y]_{\kk}+\alpha(x,y) + J[j(x)+\beta(x), y] + J [x, j(y) + \beta(y) ]\\
& =  & [x,y]_{\kk}+\alpha(x,y) + J([j(x),y]_{\kk} - \alpha(j(x), y) -\pi(y) \beta(x)) \\
& & +J([x,j(y)]_{\kk}+\alpha(x,j(y)) +\pi(x) \beta(y) )\\
& =  & [x,y]_{\kk} + \alpha(x,y) + j[j(x),y]_{\kk} + \beta([j(x),y]_{\kk}) + J\alpha(j(x),y)- J\pi(y) \beta(x)) \\
& & + j[x,j(y)]_{\kk}+ \beta([x,jy]_{\kk}) + J\alpha(x,j(y))+ J \pi(x) \beta(y)).
\end{array}
$$
\end{itemize}
Comparing both expressions  one can see that the $\kk$-component of the equality above must satisfy
\begin{equation}\label{Jink}
[j(x), j(y)]_{\kk} = [x,y]_{\kk}+ j[j(x),y]_{\kk} + j[x,j(y)]_{\kk}
\end{equation}
for an almost $j\in \End(\kk)$ such that $j^2=-1$.

While the $\hh$-component  involves  all the elements: $[\,,\,]_{\hh}, \alpha, \beta, j$ and $J$:
$$
\alpha(j(x),j(y))-\alpha(x,y)-J\alpha(x,j(y)) -J\alpha(j(x),y)=$$
\begin{equation}\label{Jinh}
\begin{array}{rcl}
& = &  \beta([j(x),y]_{\kk}) + \beta([x,jy]_{\kk}) - [\beta(x), \beta(y)]_{\hh}+ \\
& & 
+ \pi(j(y)) \beta(x) - \pi(j(x))\beta (y)  - J\pi(y) \beta(x))+ J \pi(x) \beta(y). 
\end{array}
\end{equation}

Conversely one has
\begin{prop} \label{coh}
Let $(\hh, J)$ denote a  Lie algebra equipped with a complex structure $J$. Let $\ggo =\kk \oplus \hh$ denote a Lie algebra such that $\hh$ is an ideal of $\ggo$ and $\kk$ is a linear subspace and  let $(\pi, \alpha)$ be the elements  arising from the exact sequence
$$0 \longrightarrow \hh \longrightarrow \ggo \longrightarrow \ggo/\hh \longrightarrow 0$$
as in Proposition \ref{p2}.

Let $j: \kk \to \kk$ denote an almost complex structure on $\kk$ and define $\tilde{J}: \ggo \to \ggo$ by
$$
\begin{array}{rclll}
\tilde{J}x & = & Jx & \mbox{ for } x\in \hh \\
\tilde{J}x & = & jx + \beta(x) & \mbox{ for } x\in \kk
\end{array}
$$
where $\beta\in Hom(\kk,\hh)$. Then $\tilde{J}$ defines a complex structure on $\ggo$ if and only if (\ref{jbeta}), (\ref{jkh}), (\ref{Jink}),  (\ref{Jinh})  hold.
\end{prop}

Note that for different reasons the subspace $\kk$ above could  not be $J$-invariant. For instance in presence of a symplectic structure it could be necessary to take it isotropic but not complex. 

Assume now the subspace $\kk$ is $J$-invariant. Then the starting point for the construction is the above one with $\beta=0$.  Let $J_1: \kk \to \kk$ an almost complex structure and $J_2: \hh\to \hh$ also an almost complex structure. Let $\ggo=\kk \oplus \hh$ as in Proposition \ref{lie}. The linear map $J_{\pm}=(J_1, \pm J_2)$ defines an almost complex structure on $\ggo$. 

Thus the Nijenhuis tensor on $\kk$ and $\hh$ gives that
\begin{itemize}
\item $N_{J_{\pm}}(x,y)=N_{J_1}(x,y)$ for all $x,y\in \kk$;
\item $N_{J_{\pm}}(x,y)=N_{J_2}(x,y)$ for all $x,y\in \hh.$
\end{itemize}

 
 




The proof of the following corollary  follows from the situation above in the case $\beta=0$.

\begin{cor} \label{p1} Let $\ggo=\kk \oplus \hh$ be a Lie algebra as in (\ref{lie}) attached to $(\pi, \alpha)$ and let $J_1$ denote a complex structure on $\kk$ and $J_2$ a complex on $\hh$. The almost complex structure on $\ggo$ given by 
$$J_{\pm}(x,y):= (J_1x, \pm J_2)$$ is integrable if and only if the following conditions hold

\vskip 3pt
{\rm (i)} $\varepsilon [\pi(J_1 x), J_2] y+[\pi(x), J_2] J_2 y= 0$ \quad for $x\in \kk, y \in \hh$;

\vskip 3pt

{\rm (ii)} $\alpha(J_1 x , J_1 y )- \alpha (x,y) + \varepsilon J_2 (\alpha(J_1 x, y) +  \alpha(x,J_1 y))=0$ \quad for $x, y\in \kk$,

where $\varepsilon =1$ for $J_+$ and $\varepsilon =-1$ for $J_-$.
\end{cor}

\begin{defn} \label{d1} Let  $J_1$ be a complex structure on $\kk$ and $J_2$ be  a complex structure on  a vector space $V$. Assume $\kk$ acts  on $V$ via $\pi$. We shall say that the action is {\em holomorphic} if $[\pi(x), J_2]=0$ for all $x\in \kk$.

Let $B: \kk \times \kk \to V$ be a bilinear map. We say that $B$ is {\em compatible} with $J_1$ if $B(J_1 x, J_1 y)=B(x,y)$.
\end{defn}

\begin{cor} Let $\ggo=\kk \oplus \hh$ be a Lie algebra as in (\ref{lie}) attached to $(\pi, \alpha)$. Let $J_1$ denote a complex structure on $\kk$ and $J_2$ a  complex structure on $\hh$. Assume that the action of $\kk$ into $\hh$ is holomorphic and $\alpha$ is compatible with $J_1$. Then the almost complex structure  $J_{\pm}(x,y):= (J_1 x, \pm J_2 y)$ is integrable  on $\ggo$.
\end{cor}

At the Lie group level one has the next result. See \cite{Roge} for more details. 

\begin{lem} \cite{Roge} Let $(\ggo, J)$ be a  Lie-algebra with complex structure. Let $\hh\subset \ggo$ be an ideal of $\ggo$ such is complex.  Let $G$ and $H$ denote the  associated simply connected Lie-groups endowed with the left-invariant complex structures induced by $J$ and assume thta $H$ is closed in $G$.  Then there is a holomorphic fibration  $\rho:G \to G/H$ with fiber $H$.
\end{lem}

\vskip 5pt

{\bf A totally real ideal.} Now we study  complex structures $J$ on a Lie algebra of the form $\ggo=\kk \oplus \hh$ where $\hh$ is an ideal in $\ggo$  and such that $J\kk = \hh$. 

Examples of this can be constructed from 1-cocycles as we show below. 
 Let $\kk$ be a Lie algebra  and let $\hh$ be a 2-step nilpotent Lie algebra of the same dimension as $\kk$. Let $\pi$ be a representation of $\kk$ into $\hh$ by derivations and let $j \in Z^1(\kk, \pi)$ be of  maximal rank.

Let $\alpha_j:\kk \times \kk \to \zz(\hh)$ be the skew symmetric bilinear map given by
$$\alpha_j(x,y) = [jx, jy]_{\hh}$$

This is a 2-cocycle, in fact using that $j$ is a 1-cocycle and $\pi$ acts by derivations, for  $x,y,z\in \kk$ one has
$$\begin{array}{rcl}
d\alpha_j(x,y,z) & = & [j[x,y]_{\kk}, jz]_{\hh} +  [j[y,z]_{\kk}, jx]_{\hh} + [j[z,x]_{\kk},jy]_{\hh}+ \\
& & -\pi(x) [jy,jz]_{\hh} - \pi(y) [jz,jx]_{\hh} -\pi(z) [jx,jy]_{\hh}\\
& = &   [ \pi(x) jy - \pi(y) jx, jz]_{\hh} + [\pi(y) jz -\pi(z) jy, jx]_{\hh}+\\
& & [jy, \pi(z) jx -\pi(x) jz]_{\hh}  -[\pi(x) jy,jz]_{\hh} - [jy,\pi(x)jz]_{\hh} +\\
 && -[ \pi(y) jz,jx]_{\hh}
 - [jz,\pi(y) jx]_{\hh} -[\pi(z) jx,jy]_{\hh} - [jx, \pi(z) jy]_{\hh}\\
& = & 0.
\end{array}
$$
Thus Proposition \ref{lie} says that $\ggo=\kk \oplus \hh$ attached to $(\pi,\alpha_j)$ is a Lie algebra with Lie bracket $[\cdot , \cdot]$ as in (\ref{corchete}). Moreover the  almost complex structure  $J:\ggo \to \ggo$   given by
\begin{equation}\label{complex}
J_{|_{\kk}}=j \qquad \qquad \qquad J_{|_{\hh}}=-j^{-1}
\end{equation}
is integrable. In fact, by calculating the Nijenhuis tensor $N_J(x,y)$ for $x,y\in \kk$ one gets 
\begin{equation}\label{nk}
\begin{array}{rcl}
N_J(x,y) & = & [jx,jy]_{\hh} -[x,y]_{\kk} - \alpha_j(x,y) -j\pi(x)jy -j\pi(y)jx\\
& = & [jx,jy]_{\hh} -[x,y]_{\kk} - [jx,jy]_{\hh} -j\pi(x)jy -j\pi(y)jx\\
& = & 0
\end{array}
\end{equation}
where the last equality follows from the  condition of $j$ being a 1-cocycle. These considerations prove  the following result.

\begin{prop} \label{p22} Let $\ggo=\kk \oplus \hh$ be a Lie  algebra attached to $(\pi, \alpha)$ as in Proposition \ref{lie} where $\hh$ is two-step nilpotent and $\dim \hh=\dim \kk$. Let $j\in C^1(\kk, \pi)$. Then the endomorphism $J:\ggo \to \ggo$, given by
\begin{equation}\label{j}
J(x,y)=(-j^{-1} y, jx)\qquad \qquad x,y \in \kk
\end{equation}
defines a complex structure on $\ggo$ if and only if $j$ is a 1-cocycle of maximal rank and the 2-cocycle $\alpha$ satisfies $\alpha(x,y)=[jx,jy]_{\hh}$ for all $x,y \in \kk$.
\end{prop}

The converse of the previous construction is given in the following theorem.
 
\begin{thm}\label{t1} Let $\ggo$  be a Lie algebra  with a complex structure $J$ and assume that $\ggo$ decomposes into a direct sum of vector spaces $\ggo=\kk \oplus J\kk$ where $J\kk$ is an ideal in $\ggo$.  Then  $\hh:=J\kk$ is either 2-step nilpotent or abelian and  $J$ is induced from a 1-cocycle  $j \in Z^1(\kk, \pi)$ of maximal rank , if $\kk$ is equipped with the Lie bracket of $\ggo/\hh$ and $\pi$ is a representation from $\kk$ into $\hh$ by derivations.
\end{thm}

\begin{proof} Let $\hh:=J\kk$ denote the ideal in $\ggo$. According to Proposition \ref{p2} The Lie algebra $\ggo$ is isomorphic to the Lie algebra $\kk\oplus \hh$ attached to $(\pi, \alpha)$, where $\pi(x)y=[x,y]$, for $x\in \kk$, $y\in \hh$ and $\alpha(x,y)=p_2 [x,y]$ for $x,y \in \kk$,  the linear map $p_2:\ggo \to \hh$ is the projection onto $\hh$ with respect to the decomposition $\kk \oplus \hh$. 

If $J$ is a complex structure on $\ggo$, by restricting it to $\kk$, the map $j:=J_{|_{\kk}}$ naturally induces an element $j\in C^1(\kk, \pi)$ of maximal rank, since $\ggo$ decomposes as a direct sum $\ggo=\kk \oplus J\kk$. Furthermore, since $J^2=-1$, for any $x\in \kk$ we have $- x = J^2 x = -j^{-1} jx$. Thus one can write $J$ in the form (\ref{j}).

By computing the Nijenhuis tensor  for $x,y\in \kk$ one gets that $j$ is a 1-cocycle and $\alpha(x,y)=[jx,jy]_{\hh}$ for $x,y\in \kk$. Since $Im \alpha \subset \zz(\hh)$, one obtains that the commutator of $\hh$ is contained in the center $C^1(\hh)\subset \zz(\hh)$, and this says  $\hh$ is either two-step nilpotent or abelian if $\alpha\equiv 0$.
\end{proof}

\begin{rem} \label{re1} The results here generalize those in \cite{CO}. In fact the results proved there was the following for tangent Lie algebras. Let $\ct \kk$ denote the tangent Lie algebra of a Lie algebra $\kk$. In \cite{CO} a complex structure $J$ on $\ct \kk$ such that $J\kk=\hh$ is called a totally real complex structure.  
\end{rem}

\begin{thm} \cite{CO} Let $\ct \kk$ denote the tangent Lie algebra of a Lie algebra $\kk$. The set of totally real complex structures on $\ct \kk$ is in one to one correspondence with the set of non-singular derivations of $\kk$. 

If one set of those (and therefore both) is non-empty then $\kk$ is nilpotent. 
\end{thm}

\begin{rem} Totally real complex structures on semidirect products of the form $V \rtimes_{\pi} \kk$, where $V$ is the underlying vector space to $\kk$ equipped with its canonical abelian bracket, are in correspondence with Lagrangian symplectic structures on $V \rtimes_{\pi^*} \kk$. See e.g. \cite{CLP}. 
\end{rem} 

\begin{exa} In the paragraphs we get examples of an {\em (almost) complex product structure} on a Lie algebra $\ggo$: that is a pair $(J,E)$ of an (almost) complex structure $J$ and an (almost) product structure $E$  (Remark \ref{prods}) such that $JE = - EJ$. See \cite{AS,BV}.
\end{exa}
\section{Bilinear forms, ideals and complex structures}

In this section we shall study compatibility conditions for $J$ with respect to non-degenerate either symmetric or skew-symmetric bilinear forms  on a Lie algebra $\ggo$ having a fixed ideal. 

\medspace

\subsection{Symmetric case}

A {\em metric}  on a real vector space $\vv$ is a symmetric bilinear map on $\vv$, $\la \, , \,  \ra:\vv \times \vv \to \RR$ which is non-degenerate, that is, for any non zero vector $x\in \vv$ there exists a vector $y\in\vv$ such that $\la x, y\ra \neq 0$. Otherwise  $\la \, , \, \ra$ is said degenerate.
   
   If $\ww$ is a subspace of  $(\vv, \la \, , \, \ra)$ the subspace
 $$\ww^{\perp} = \{ x \in \vv \,:\, \la x, v\ra =0 \quad \text{ for all } \quad v \in \ww\}$$
 denotes the {\it orthogonal}  subspace of $\ww$. In particular we say that $\ww$ is 
\begin{itemize}
\item {\em isotropic} if 
 $\ww \subset \ww^{\perp}$, 
\item {\em totally isotropic} if $\ww =\ww^{\perp}$  and 
\item {\em non-degenerate} if $\ww \cap \ww^{\perp}=0$. 
\end{itemize}

  In the last case, the restricted metric $\la \, , \, \ra_{\ww}:=\la \, , \, \ra_{|_{\ww \times \ww}}$, defines a isomorphism $\xi$ between $\ww$ and its dual space $\ww^{\ast}$ by $\xi(u)(v)=\la u,v\ra_{\ww}$ whenever $\ww$ is  finite dimensional. As usual, a metric of index 0 or signature $(0,n)$ on a vector space of dimension $n$ is called an {\it inner product}.
 
 \begin{exa} \label{hyp} Let  $\uu$ denote a vector space  whose dual space is  denoted by $\uu^{\ast}$. Let $\uu\oplus \uu^*$ be the direct sum as vector spaces of $\uu$ and $\uu^*$ and endow this with the hyperbolic metric $\la x_1 +\phi_1, x_2 + \phi_2 \ra = \phi_1 (x_2) + \phi_2(x_1)$ where $\phi_i \in \uu^*$, $x_i\in \uu$, for  i=1,2. Clearly  $\uu$ and $\uu^*$ are complementary totally isotropic subspaces in $\uu \oplus \uu^*$. 
 \end{exa}

\begin{defn} Let $\la \,,\, \ra$ denote a metric on a Lie algebra $\ggo$. Let $J$ denote an (almost) complex structure on $\ggo$. The pair $(J, \la\,,\ra)$ defines  
\begin{itemize}
\item an Hermitian structure  on $\ggo$ if $\la Jx, Jy\ra = \la x, y\ra$ $\forall x, y\in\ggo$,

\item an (almost) anti-Hermitian structure  on $\ggo$ if $\la Jx, Jy\ra = -\la x, y\ra$ $\forall x,y\in\ggo$.
\end{itemize}

We shall also say that the metric $\la\,,\,\ra$ is (almost)-Hermitian or (almost) anti-Hermitian.

We note that here we consider  Hermitian structures in relation to possibly definite metrics, although this notion is referred to by other authors as definite metrics.  Anti-Hermitian structures are also called Norden metrics \cite{No} or B-metrics.

\end{defn}

In the next paragraphs we  discuss different  possible constructions. 

\medskip

{\bf Case of a complex ideal.} 
Here we shall consider different metric constructions for the complex ideal $\hh$.

\smallskip

 Let $\ggo$ denote a Lie algebra with a complex structure $J$ and let $\hh\subset \ggo$ be a complex ideal on $\ggo$. Assume $\la\,,\,\ra$ is a metric on $\ggo$ and $\hh$ is non-degenerate relative to the metric. Let $\kk$ be  the orthogonal complementary subspace of $\hh$. Then as in \ref{coh}, the complex structure $J$ does not need to leave $\kk$ invariant. The pair $(J,\la\,,\,\ra)$ is Hermitian if and only if
$$\la Jx, Jy\ra=\la x,y\ra \qquad \forall x,y \in \ggo$$
which implies
\begin{enumerate}
\item $\la Jx, Jy\ra=\la x,y\ra \qquad \forall x,y \in \hh$
\item 

$\begin{array}{rcl} \label{comp1}
\la j(x)+\beta(x), j(y)+\beta(y)\ra & = & \la j(x), j(y)\ra+\la\beta(x), \beta(y)\ra\\
& = & \la x, y\ra \qquad \forall x,y \in \kk
\end{array}
$
\end{enumerate}

which says that $(J_{|_{\hh}}, \la\,,\,\ra_{\hh})$ defines an (almost) Hermitian structure on $\hh$ and 
on $\kk$ one should have:
$$\la j(x), j(y)\ra - \la x, y\ra +\la\beta(x), \beta(y)\ra=0 \qquad \mbox{for all } x,y \in \kk.$$
 
Thus in this situation $\hh$ is a complex non-degenerate ideal and $\kk$ is non-degenerate albeit $J$-invariant. 

\smallskip

Some cases are the following ones.
\begin{enumerate}
\item   Let $\kk$ and $\hh$ denote Lie algebras with corresponding Hermitian structures $(J_1, B_1)$ and $(J_2, B_2)$. Let $\ggo=\kk \oplus \hh$ be a Lie algebra as in Proposition\ref{lie}  attached to $(\pi, \alpha)$, then $(J_{\pm}, B_1 + B_2)$ defines an almost Hermitian structure on $\ggo$.

\item Let $(\kk, J_1)$ and $(\hh, J_2)$ denote Lie algebras with corresponding (almost) complex structures. Assume $\dim \hh=\dim \kk$ and let $t:\kk \to \hh$ be an isomorphism.  Let $B$ be a metric on $\kk$ and consider the metric on $\ggo$ given by
\begin{equation}\label{m2}
\la (x_1, ty_1), (x_2, t y_2)\ra = B(x_1, y_2)+B(x_2, y_1).
\end{equation}
Then the pair $(J:=(J_1, \pm J_2), \la\,,\,\ra)$ defines an (almost) Hermitian structure on $\ggo$ if and only if $J_2 t= \pm t J_1$.

Notice that in this situation both $\kk$ and $\hh$ are isotropic subspaces and the metric has signature $(n,n)$ where $n=\dim \kk$. 
\end{enumerate}

\begin{cor} Let $\kk$ denote a Lie algebra with an (almost) Hermitian structure $(J_1, B)$. Let $\ggo=\kk \oplus \hh$ be the extended Lie algebra such that $\hh$ is an ideal of $\ggo$. Assume $\dim \hh=\dim \kk$ and let $t:\kk \to \hh$ be a linear isomorphism. Define an almost structure $J_2$ on $\hh$ by 
$$J_2 = t J_1 t^{-1}$$
Then the metric on $\ggo$ given by
$\la (x_1, ty_1), (x_2, ty_2)\ra = B(x_1, y_2) + B(x_2, y_1)$
gives rise to an (almost) Hermitian structure for $J=(J_1, J_2)$.
\end{cor}

\smallskip

{\bf Case of a totally real ideal.}
Here we shall consider different metric constructions for the totally real  ideal $\hh$.

\smallskip

Let  $\ggo$ be a Lie algebra which splits into a direct sum of vector spaces $\ggo=\kk \oplus \hh$  and which admit  an almost complex structure $J$ such that $J\kk=\hh$: i.e. there is a linear isomorphism $j: \kk \to \hh$ (see (\ref{lie})). 
$$J(x,y)=(-j^{-1}y, jx)\qquad \mbox{for all} x\in \kk, y\in \hh.$$
Let $B$ denote a metric on $\kk$. Consider the following metric $\la \, , \, \ra$ as extension of $B$ to $\ggo$:
\begin{enumerate}
\item  
$$\la (x_1, j y_1), (x_2, j y_2)\ra = B(x_1, x_2) + B (y_1, y_2)\qquad \mbox{ for all } x_i, y_i \in \kk, \, i=1,2.$$
This metric  is Hermitian:
$$\begin{array}{rcl}
\la J(x_1, j y_1), J(x_2, j y_2)\ra & = & \la (-y_1, j x_1), (-y_2, j x_2)\ra\\
 & = & B (y_1,  y_2) + B (x_1,x_2)
 \end{array}
 $$
Clearly $\la \, , \, \ra$ restricts to both $\kk$ and $\hh$ as a metric and so that $\kk \perp \hh$. However the geometry that $\la \, \, \ra$ determines on $\kk$ is different from that one on $(\kk, B)$.  
 
 Let us denote by $\nabla$ the Levi  Civita connection corresponding to $\la \, , \, \ra$.
For $x,y, z \in \kk$ the following formulas hold
$$2 \la \nabla_x y, jz \ra = \la \alpha(x,y), z\ra\qquad  2\la \nabla_{jx} jy, z\ra=\la \pi(z) jy, jx\ra - \la \pi(z) jx, jy\ra$$
showing that $\kk$ and $\hh$ are not necessarily totally geodesic subspaces. 

\item $$\la (x_1, j y_1), (x_2, j y_2)\ra = B(x_1, y_2) + B (x_2, y_1)\qquad \mbox{ for all } x_i, y_i \in \kk, \, i=1,2.$$
This metric  is anti-Hermitian. Both spaces $\kk$ and $\hh$ are isotropic. 
\end{enumerate}

\begin{prop}\label{p3} Let $(\ggo, J)$ be a Lie algebra equipped with an almost complex structure $J$ and assume $\ggo$ splits into a direct sum of vector spaces $\ggo=\kk \oplus \hh$ such that $J\kk =\hh$. Then $\ggo$ always admits an Hermitian and an anti-Hermitian metric for $J$.
\end{prop}
\begin{proof} Let $B$ denote an inner product on $\kk$ and denote by $j$ the restriction of $J$ to $\kk$, $j:=J_{|_{\kk}}$. Since $J\kk=\hh$,  the almost complex structure $J$ induces a linear morphism $j:\kk \to \hh$ which is non singular and since $J^2=-1$, it is easy to see that $J$ is related to $j$ by the formula  $J(x,y)=(-j^{-1} y, jx)$.

 Let $\la \, , \, \ra$ be the metric on $\ggo$ given by
$$\la (x_1, y_1), (x_2, y_2)\ra = -B(x_1, j^{-1} y_2) - B(x_2, j^{-1} y_1).$$
The map $\la \, , \, \ra$ is bilinear  symmetric and non-degenerate. Moreover it is compatible with $J$, as it follows from
$$\begin{array}{rcl}
\la J(x_1, y_1), J(x_2, y_2)\ra & = &\la (-j^{-1}y_1, j x_1), (-j^{-1}y_2, j x_2)\ra \\
& = & -B(-j^{-1}y_1,  x_2) - B(-j^{-1}y_2, x_1)\\
& = & \la (x_1, y_1), (x_2, y_2)\ra.
\end{array}$$
Note that the subspaces $\kk$ and $\hh$ are isotropic and of maximal dimension hence they are totally isotropic.

For the anti-Hermitian structure  the metric on $\ggo$ is defined as 
$$\la (x_1, y_1), (x_2, y_2)\ra = -B(x_1, j^{-1} y_2) + B(x_2, j^{-1} y_1)$$
and the proof follows along  the same lines  of the preceding case.
\end{proof}


\smallskip

{\bf A remark on  SKT structures}
Let $(M, J, g)$ be a Hermitian manifold. If the torsion 3-form $c$  of the Bismut
connection is $d$-closed, then the Hermitian metric $g$ on a complex manifold $(M, J)$ is called
strong K\"ahler with torsion (shortly SKT), where
$$c(x,y,z)=g(x, T^B(y,z))$$
being $T^B$ the torsion of the Bismut connection $\nabla^B$ characterised as  the unique connection on the Hermitian manifold $(M,J,g)$ such that $\nabla^BJ=0$, \, $\nabla^Bg=0$.

Let $\sso(\ggo,g)$ denote the Lie algebra of skew symmetric maps for $g$. See \cite{EF} for the proof of the next result and for more details and references on SKT structures. 

\begin{prop} \cite{EF} Let $(\ggo, J, g)$ be a Hermitian Lie algebra and let $\pi: \ggo \to \sso(\ggo,g)$ be a representation such that $\pi$ is holomorphic. Take $\ggo \ltimes_{\pi}\hh$ where $\hh$ is the vector space underlying $\ggo$ with the trivial bracket.  Then the Hermitian structure
$(\tilde{J}, \tilde{g})$ given by $\tilde{g}=g+g$ the product metric on $\ggo \oplus \hh$ and $\tilde{J}(x,y)=(Jx,Jy)$ is SKT if and only if $(J, g)$ is SKT on $\ggo$.
\end{prop}

In remark 3.2 the authors exemplify the result for the adjoint representation. They say that the conditions of the Proposition hold in this situation if and only if on the Hermitian Lie algebra $(\ggo,J,g)$ the complex structure $J$ is bi-invariant and the inner product $g$ is ad-invariant. However as proved in \cite{ABO} this is possible only for an abelian Lie algebra $\ggo$.

\subsection{Skew-symmetric case} Here we shall study symplectic structures on Lie algebras $\ggo$ with an ideal $\hh$.

\vskip 3pt

A {\em symplectic structure}  on a Lie algebra $\ggo$ is a skew-symmetric non-degenerate bilinear form on $\ggo$, $\omega:\ggo \times \ggo \to \ggo$ satisfying the {\em closeness} condition:
$$\omega([x,y], z)+\omega([y,z],x)+\omega([z,x],y)=0\qquad \mbox{ for all }x,y,z\in \ggo.$$

The pair $(\ggo, \omega)$ is sometimes called a {\em symplectic Lie algebra.}
   
   If $\ww$ is a subspace of  $(\ggo, \omega)$ the subspace
 $$\ww^{\perp_{\omega}} = \{ y \in \ggo \,:\, \omega( x, y) =0 \quad \text{ for all } \quad x \in \ww\}$$
 denotes the {\it symplectic-orthogonal}  subspace of $\ww$.
 
  In particular we say that $\ww$ is
  
  \begin{itemize}
\item    {\em isotropic} if  $\ww \subset \ww^{\perp_{\omega}}$, 
 \item {\em Lagrangian} if $\ww =\ww^{\perp_{\omega}}$,
 \item  {\em symplectic} if $\ww \cap \ww^{\perp_{\omega}}=0$. 
\end{itemize}

\begin{prop}
Let $\hh$ be an isotropic ideal on a symplectic Lie algebra $(\ggo, \omega)$ then $\hh$ is abelian.
\end{prop}
\begin{proof} The closedness condition for $x,y\in \hh, z\in \ggo$ 
$$0=\omega([x,y],z)+\omega([y,z],x)+\omega([z,x],y)$$
and since $[y,z]\in \hh\subset \hh^{\perp_{\omega}}$ one gets $\omega([y,z],x)=0$. 

Analogously $\omega([z,x],y)=0$. Therefore
$$0=\omega([x,y],z)\qquad \mbox{ for all }z\in \ggo$$
and since $\omega$ is non-degenerate, one has $[x,y]=0$ for all $x,y\in \hh$.
\end{proof}

Particular examples arise in the next context. Let $J$ be a complex structure on a symplectic Lie algebra $(\ggo, \omega)$. The pair $(\omega, J)$ is called 
\begin{itemize}
\item a pseudo-K\"ahler structure on $\ggo$ if $\omega(Jx,Jy)=\omega(x,y)$ for all $x,y\in \ggo$.
\item a complex-symplectic structure on $\ggo$ if
$\omega(Jx,Jy)=-\omega(x,y)$ for all $x,y\in \ggo$. See \cite{COP} for this definition.
\end{itemize}

\begin{rem} Notice that there is a one-to-one correspondence between non-degenerate skew-symmetric bilinear maps $\omega$ compatible with an almost complex structure $J$ and metrics $\la\,,\,\ra$ compatible with an almost complex structure $J$ due to 
$$\omega(x,y)=\la x, Jy\ra \qquad \forall x,y \in \ggo.$$
Assume $J$ is integrable. Then if $\nabla$ denotes the Levi-Civita connection associated to $\la\,,\,\ra$ then $\nabla J=0$ is equivalent to $d\omega=0$.

Thus for a complex structure $J$ one gets a one-to-one correspondence between

\smallskip

$$\begin{array}{lll}
\left\{
{\begin{array}{l}
\mbox{Hermitian structures}\\ \mbox{ such that }\nabla J=0 
\end{array}}\right\}
  &   \longleftrightarrow  & 
	\left\{ 
	\begin{array}{l}
\mbox{ pseudo-K\"ahler structures for } J \end{array}\right\}\\
 & &  \\ \\
\left\{
{\begin{array}{l}
\mbox{anti-Hermitian structures}\\ \mbox{ such that } \nabla J=0
\end{array}} \right\}  & \longleftrightarrow  & \left\{\mbox{ complex-symplectic structures  for }J\right\}
\\
 & &  \\
\end{array}
$$
\end{rem}

So this equivalence and Proposition \ref{p3} gives the next result.

\begin{cor} Let $(\ggo, J)$ be a Lie algebra equipped with an almost complex structure $J$ and assume $\ggo$ splits into a direct sum of vector spaces $\ggo=\kk \oplus \hh$ such that $J\kk =\hh$. Then $\ggo$ always admits a non-degenerate skew-symmetric bilinear map which is compatible or anti-compatible with $J$.
\end{cor}

\begin{exa} Let $\kk$ be a Lie algebra endowed with  an (almost) Hermitian structure $(B,J_1)$ and let $\hh$ denote a Lie algebra endowed with an almost Hermitian structure $J_2$. Let $\ggo=\kk \oplus \hh$ denote the direct sum of vector spaces with complex structure $J:=(J_1, \pm J_2)$ as in Remark \ref{p1}. Assume that $\dim \kk=\dim \hh$ and let $t:\kk \to \hh$ be an isomorphism. If $tJ_1=\pm J_2 t$, then the bilinear map $\Omega:\ggo \times \ggo \to \RR$ given by
$$\Omega((x_1, ty_1), (x_2, ty_2))= B(x_1, J_1 y_2)- B(x_2, J_1 y_1)$$
  is skew-symmetric and compatible with $J$.
\end{exa}

\begin{prop}
Let $J$ denote a complex structure on a symplectic Lie algebra $\ggo$. Then if $\hh$ is an isotropic ideal, $J\hh$ is a Lie subalgebra of $\ggo$. 

Moreover if $\hh$ is totally real, then $\ggo=J\hh \ltimes \hh$.
\end{prop}
\begin{proof} The previous proposition says that $\hh$ must be abelian. Now the integrability condition for $J$ gives
$$[Jx,Jy]=J[Jx,y]+J[x,Jy] \qquad \mbox{ for all } x,y \in \hh$$
and $[Jx,y]=\pi(Jx) y\in \hh$ so as $[x,Jy]=-\pi(Jy) x \in \hh$ which says that $\kk=J\hh$ is a Lie subalgebra of $\ggo$. 

If $\hh$ is totally real then $\hh \cap J\hh=\{0\}$ and $\ggo=J\hh \ltimes \hh$.
\end{proof}

\begin{rem}
In the situation of the preceding proposition, if $\hh$ is Lagrangian and the pair $(\omega,j)$ is either pseudo-K\"ahler or complex-symplectic, then $J\hh$ is a Lagrangian subalgebra.
\end{rem}

Let $\kk$ denote a Lie algebra and let $\gamma: \kk \to \End(\kk)$ denote a linear map. Then we say that $\gamma$ is a connection which is
\begin{itemize}
\item {\em torsion-free } if $\gamma(x) y - \gamma(y) x=[x,y]_{\kk}$
\item {\em flat} if $\gamma([x,y])=\gamma(x)\gamma(y)-\gamma(y)\gamma(x)$
 for all $x,y\in \kk$, that is $\gamma: \kk \to \End(\kk)$ is a representation of Lie algebras. 
 \item Given a symplectic structure on $\kk$ the connection $\gamma$ on $\ggo$ is said to be {\em symplectic } if 
 $$\omega(\gamma(x) y, z)+ \omega(y, \gamma(x) z)=0$$
 \end{itemize}
 
 Lie algebras $\ggo$ endowed with a symplectic structure and a torsion-free flat symplectic connection give rise to hypersymplectic structures on $\ggo \ltimes V$ where $V$ is the underlying vector space to $\ggo$, as proved in \cite{COP}. (See proofs and definitions there).
 
\begin{thm}\cite{COP} Let $(\ggo, \omega)$ be a symplectic Lie algebra with a torsion free, flat symplectic connection $\gamma$ on the underlying vector space $V$ of the Lie algebra. Then the associated space $\ggo \ltimes V$ admits an hypersymplectic structure such that the Levi-Civita connection of the associated neutral metric is flat and symplectic with respect to each of the three given symplectic structures.
\end{thm}

\smallskip

{\bf Generalized complex structures} Let $\kk$ be a Lie algebra, and let $\kk^*$ denote its dual space. The cotangent Lie algebra $\ct^*\kk$ is the semidirect product of $\kk$ and $\kk^*$ via the coadjoint action $\ad^*:\kk \to \End(\kk^*)$, which is given by
$$ \ad^*(x) \cdot \varphi = - \varphi \circ \ad(x) \qquad \mbox{ for } x\in \kk, \varphi \in \kk^*.$$

The canonical neutral metric on $\ct^*\kk$ is that one already defined in Example (\ref{hyp}), also called hyperbolic metric:
$$ \la (x_1, \varphi_1), (x_2, \varphi_2)\ra=  \varphi_1 (x_2) + \varphi_2 (x_1).$$

According to \cite{ABDF} a {\em generalized complex structure} on a Lie algebra $\kk$ is an Hermitian structure $(J,\la \, , \, \ra)$ on its cotangent Lie algebra $\ct^*\kk$, where $\la \, , \, \ra$ is the canonical neutral metric. 

One can see that the dimension of $\kk$ must be even.
Assume  the dimension of $\kk$ is $2n$, $\dim \kk=2n$. If we choose on $\ct^* \kk$ a basis adapted to the splitting $\kk \oplus \kk^*$, the matrix of $J$  has the following form
$$ J= \left( \begin{matrix} j_1 & j_2 \\ j_3 & j_4
\end{matrix} \right)
$$
where $j_i$ are certain $2n\times 2n$ matrices for  $i=1,2,3,4$. One says that $J$ is of {\em type k} if $\rank j_2 = 2(n-k)$. When $j_2=j_3 =0$ the generalized complex structure $J$ is said to be of {\em complex type}; if $j_1= j_4=0$ the generalized complex structure is called of {\em symplectic type}. 

If the Lie algebra $\kk$ itself is endowed with a complex structure, then $\kk$ has a generalized complex structure. In fact, if $J_1:\kk \to \kk$ is a complex structure on $\kk$, extending it to $\kk^*$ as $J_2(\varphi) = -\varphi \circ J_1$, the almost complex structure on $\ct^*\kk$ given by $J(x,\varphi)=(J_1 x, J_2 \varphi)$ is integrable and compatible with the canonical neutral metric. Notice that this is a particular case of (\ref{p1}).

Applying results of the previous section we are able to  characterize generalized complex structures of symplectic type. In fact, let  $J$ denote an Hermitian structure on $\ct^*\kk$ of symplectic type, then its restriction  to $\kk$ induces a linear morphism $j:=J_{|_{\kk}}:\kk \to \kk^*$. The integrability of $J$ says that $j$ is a 1-cocycle of $(\kk, \ad^*$).

\begin{cor} Any generalized complex structure of symplectic type on a even dimensional Lie algebra $\kk$ is determined by a 1-cocycle of $(\kk, \ad^*)$.
\end{cor}

Conversely (\ref{p22}) and (\ref{p3}) imply the following result.

\begin{cor} Let $\kk$ be a even dimensional Lie algebra and let $j:\kk \to \End(\kk^*)$ denote a 1-cocycle of $(\kk, \ad^*)$. Then $j$ induces a generalized complex structure of symplectic type on $\ct^*\kk$.
\end{cor}
\begin{proof} The almost complex structure on $\ct^*\kk$ given by $J(x,y):=(-j^{-1}y, jx)$ for $x\in \kk$ $y\in \kk^*$ is integrable (see \ref{p22}). Since Proposition \ref{p3} applies for almost complex structures, the result follows at once.
\end{proof}

\section{Examples of complex structures on six dimensional Lie algebras}

Our goal in this section is to apply Proposition \ref{p22} to construct complex structures  on six dimensional Lie algebras $\ggo=\kk \oplus \hh$ where  $\hh$ and  $\kk$ have the same dimension and  $\hh$ is a totally real ideal of $\ggo$.
 
 Recall first the classification of real  three dimensional Lie algebra $\kk=span\{e_1,e_2,e_3\}$, which are listed below  (see \cite{Mi} for instance):
$$
\begin{array}{ll}
\RR^3: & \mbox{ with trivial Lie bracket }\\
\hh_1: & [e_1,e_2]=e_3.\\
\rr_{3} & [e_1,e_2]=e_2,\quad [e_1,e_3]=e_2+e_3.\\
\rr_{3, \lambda}: & [e_1,e_2]=e_2,\quad [e_1,e_3]=\lambda e_3.\\
\rr_{3,\delta}: & [e_1,e_2]=e_2+\delta e_3\quad [e_1,e_3]=-\delta e_2+e_3.\\
\sso(3): & [e_1,e_2]=e_3,\quad [e_1,e_3]=-e_2, \qquad [e_2,e_3]=e_1.\\
\ssl(2): &  [e_1,e_2]=2, e_2\quad [e_1,e_3]=2 e_3,\quad [e_2,e_3]=e_1.
\end{array}
$$

\

{\bf  The construction.} Let $\kk$ be a three dimensional Lie algebra and let $\hh=\RR^3$ 
Let $\pi$ denote a representation from $\kk$ into $\hh$ by derivations. We search for $j\in Z^1(\kk,\pi)$ of maximal rank, that is $j:\kk \to \RR^3$ is an invertible linear operator satisfying  the linear equation
\begin{equation}\label{equ1}
0 =\pi(x) j(y)-\pi(y) j(x) - j([x,y]_{\kk})\qquad \mbox{ for all }x,y \in \kk.
\end{equation}
 
 This gives an integrable almost complex structure $J$ on $\ggo=\kk \ltimes \RR^3$, such thar both $\kk$ and $\hh$ as subspaces of $(\ggo, J)$ are totally real. 
The complex structure $J$ on $\ggo$ is defined by $J_{|_{\kk}}=j$, $J_{|_{\hh}}=j^{-1}$. So if $\mathcal B$ is a basis of $\kk$ and $\mathcal B'$ is a basis of  $\hh$, then $\mathcal B \cup \mathcal B'$ is a basis of $\ggo$ with respect to which the complex structure $J$ has the form 
 $$\left( \begin{matrix}0 & -j^{-1} \\
j & 0 \end{matrix} \right),
$$ where $j$ is a $ 3\times 3$ invertible real matrix. We shall write it with coefficients $j_{uv}$ as follows:
$$j= \left( \begin{matrix}
j_{41} & j_{42} & j_{43} \\
j_{51} & j_{52} & j_{53} \\
j_{61} & j_{62} & j_{63}
\end{matrix} \right)\qquad \quad \det j\neq 0.
$$

Next  we evaluate the equation (\ref{equ1}). If $\{e_1,e_2,e_3\}$ denotes a basis of $\kk$ and $\{e_4,e_5,e_6\}$ denotes a basis of $\hh$, Equation (\ref{equ1}) becomes

$$
\begin{array}{rcl}
0 & = & \pi(e_i) j(e_k) - \pi(e_k) j(e_i) - j [e_i,e_k]_{\kk}\\
& = &  \sum_{s=4}^6 j_{sk} \pi(e_i) e_s - \sum_{s=4}^6 j_{si} \pi(e_k) e_s
- j (\sum_{l=1}^3 C_{ik}^l) e_l
\end{array}
$$
which is a linear system on the coefficients $j_{uv}$, for a fixed representation by derivations $\pi$ of $\kk$ into $\RR^3$ and for $\{C_{ik}^l\}$ being the structure coefficients for $\kk$.  We shall choose a representation $\pi$ such  that $\Ima \pi=1$ to do explicit computations.

For $\hh =\RR^3$,  any representation $\pi:\kk \to \End(\RR^3)$  such that $\dim\Ima \pi=1$ is determined by a linear map $t$. Thus there is a basis of $\RR^3$ in which the matrix of $t$ is of one and only one of the following types

\begin{equation}\label{typet1} i)\, \left( \begin{matrix}
\eta & 0 & 0 \\
0 & \nu & 0 \\
0 & 0 & \mu
\end{matrix} \right)\qquad ii)\, \left( \begin{matrix}
\eta & 0 & 0 \\
0 & \nu & -\mu  \\
0 & \mu & \nu
\end{matrix} \right)\qquad iii)\,
\left( \begin{matrix}
\eta & 0 & 0 \\
0 & \nu & 1 \\
0 & 0 & \nu
\end{matrix} \right)\qquad iv)\, 
\left( \begin{matrix}
\eta & 1 & 0 \\
0 & \eta & 1 \\
0 & 0 & \eta
\end{matrix} \right).
\end{equation}

\vspace{.3cm}

For $\hh =\hh_1$,  any representation $\pi:\kk \to \Der(\hh_1)$  such that $\dim\Ima \pi=1$ is determined by a derivation $t$     of $\hh_1$ . Thus there is a basis of $\hh_1$ in which the matrix of $t$ is of one and only one of the following types (see for instance \cite{ABDO}):

\begin{equation}\label{typet2} i)\, \left( \begin{matrix}
\eta & \nu & 0 \\
\mu & -\eta & 0 \\
0 & 0 & 0
\end{matrix} \right)\qquad \qquad ii)\, \left( \begin{matrix}
\eta & \nu & 0 \\
\mu & 1-\eta & 0  \\
0 & 0 & 1
\end{matrix} \right).
\end{equation}

In any case if  we assume that $\pi(e_1) = \varepsilon_1 t$, $\pi(e_2)=\varepsilon_2 t$, $\pi(e_3)=\varepsilon_3 t$, the condition of $\pi$ being a representation says
$$\phi([x,y]) =[\pi(x), \pi(y)]= 0$$
where the last equality holds due to $\Ima \pi=span\{t\}$.  
Hence $\pi(x)=0$ for every $x\in C^1(\kk)$.

This explanation gives the proof of the following Lemma. See for instance \cite{Va} for representations of $\ssl(2)$. 

\begin{lem} The simple Lie algebras $\ssl(2)$ and $\sso(3)$ do not admit any representation $\pi:\kk\to\End(V)$ such that $\dim \Ima \pi=1$. 
\end{lem}

For the solvable Lie algebras $\hh_1, \rr_3, \rr_{3,\lambda}, \rr_{3,\delta}$ one has $\pi(e_j)=0$ for $e_j\in C^1(\kk)$. Since in all these cases $e_1\in\kk-C^1(\kk)$ we shall assume $\pi(e_1)=t$ and we take $\pi(e_2)=\varepsilon_1 t, \pi(e_3)=\varepsilon_2 t$, where $\varepsilon_i$ could be zero as explained above.

The computations for the proof of the next theorems can be done with help of a software. 

\begin{thm} Let $\kk$ be a three dimensional solvable Lie algebra. Let $\pi:\kk \to \End(\RR^3)$ be a representation  such that $\Ima \pi=span\{t\}$ where $t$ is as in (\ref{typet1}). Then the semidirect product Lie algebra $\ggo=\kk\oplus_{\pi} \RR^3$ admits a  complex structure $J$ such that $J\kk=\RR^3$ in the cases exposed in  Table 1.
\end{thm}

\begin{thm} Let $\kk$ be a three dimensional solvable Lie algebra. Let $\pi:\kk \to \End(\hh_1)$ be a representation acting by derivations on the Heisenberg Lie algebra $\hh_1$  with $\dim \Ima \pi= 1$ as in (\ref{typet2}).  

Then there exists a Lie algebra  with a  complex structure $(\ggo, J)$ such that $\hh_1$ is a totally real ideal of the   extended semidirect product Lie algebra $\ggo=\kk\oplus_{\pi} \hh_1$
as in  Proposition \ref{p22} in the cases exposed in Table 2. 

\end{thm}


\begin{center} Table 1 \end{center}

\smallskip

\begin{tabular}{|c|c|c|c|}
\hline
\quad $\kk$ \quad & {\small Type} $t$ & {\small Existence results } & {\small \qquad Parameters of $t$  for the existence \qquad } \\ \hline
$\hh_1$ & {\small (i)} & {\small no} &  \\ \hline
 & {\small (ii)} & {\small no} &   \\ \hline
 & {\small (iii)} & {\small yes} & {\small $\nu=0$}\\ \hline
  & {\small (iv)} & {\small yes} & {\small $\eta=0$} \\ \hline
$\rr_3$ & {\small (i)} & {\small no } & \\ \hline  
& {\small (ii)} & {\small no } & \\ \hline  
& {\small (iii)} & {\small yes } & {\small $\nu=1$} \\ \hline  
& {\small (iv)} & {\small yes } & {\small $\eta=1$} \\ \hline 
$\rr_{3,\lambda}$ & {\small (i)} & {\small yes } & {\small $\eta=\nu=0, \mu=1, \varepsilon=0$}   \\ 
$\lambda=0$ & & & {\small $\eta=\mu, \nu=1, \varepsilon=0$} \\
& & & {\small $\eta=1, \nu=\mu=0, \varepsilon=0$}\\
& & & {\small $\nu=1, \mu=0, \varepsilon=0$}\\
& & & {\small $\nu=0, \mu=1, \varepsilon=0$} \\
& & & {\small $\eta=1, \nu=0, \varepsilon=0$} \\
& & & {\small $\eta=1, \mu=0, \varepsilon=0$} \\
& & & {\small $\eta=0, \mu=1, \varepsilon=0$}\\
& & & {\small $\eta=0, \nu=1, \varepsilon=0$} \\ 
&  & & {\small $\eta=, \mu=\nu=1\varepsilon=0$} \\
& & & {\small $\eta=1, \nu=0, \mu=1,\varepsilon=0$} \\
& & & {\small $\eta=1, \nu=1, \mu=0, \varepsilon=0$} \\ 
\hline
	\end{tabular}

\begin{tabular}{|c|c|c|c|}
\hline
\quad $\kk$ \quad & {\small Type} $t$ & {\small Existence results } & {\small \qquad Parameters of $t$  for the existence \qquad } \\ \hline
 & {\small (ii)} & {\small yes } &  {\small $\eta=1, \nu=\mu=\varepsilon=0$}\\
& & & {\small $\nu=1, \eta=\mu=\varepsilon=0$} \\
\hline
 & {\small (iii)} & {\small yes } & {\small $\eta=1, \nu=\varepsilon=0$}\\
 &  &  & {\small $\nu=1,\eta=\varepsilon=0$}\\ 
\hline
& {\small (iv)} & {\small no } & \\
\hline
$\rr_{3,\lambda}$  & {\small (i)} & {\small yes } &  {\small $ \eta=1, \mu=\lambda$}\\ 
$\lambda\neq 0$ & & & {\small $\eta=1, \nu=\lambda$} \\
& & & {\small $\eta=\lambda, \nu=1$} \\
& & & {\small $\eta=\lambda, \mu=1$} \\
& & & {\small $\nu=1, \mu=\lambda$} \\
& & & {\small $\nu=\lambda, \mu=1$} \\
& & & {\small $\eta=1, \nu=\mu=\lambda$} \\
& & & {\small $\eta=\nu=1, \mu=\lambda$} \\
& & & {\small $\eta=\mu=1, \nu=\lambda$} \\
& & & {\small $\eta=\nu=\lambda, \mu=1$} \\
& & & {\small $\eta=\mu=\lambda, \nu=1$} \\
& & & {\small $\eta=\lambda, \nu=\mu=1$} \\
& & & {\small $\lambda=\eta=\nu=1$} \\
& & & {\small $\lambda=\eta=\mu=1$} \\
& & & {\small $\lambda=\nu=\mu=1$} \\
& & & {\small $\lambda=\eta=\nu=\mu=1$} \\ \hline
& {\small (ii)} & {\small yes } &  {\small $\lambda=\nu=1, \mu=0 $} \\
& & &  {\small $\eta=1, \nu=\lambda, \mu=0 $}\\
& & &  {\small $\eta=\lambda, \nu=1, \mu=0 $}\\
 \hline
& {\small (iii)} & {\small yes } &  {\small $\eta=1, \nu=\lambda$}\\
& & &  {\small $\eta=\lambda, \nu=1  $}\\
& & &  {\small $\lambda=\eta=\nu=1  $}\\
 \hline
& {\small (iv)} & {\small no } & \\
\hline 
$\rr_{3,\delta}$ & {\small (i)} & {\small yes} & {\small $\delta=0, \nu=0, \mu=1$ }  \\ 
  & & & {\small $\delta=0, \eta=1, \nu=0 $} \\
\hline
  & {\small (ii)} & {\small yes} & {\small $\nu=1, \mu=\delta$} \\ 
& & & {\small $\nu=1, \mu=-\delta  $} \\
\hline
& {\small (iii)} & {\small yes} & {\small $\delta=0, \eta=\nu=1$}\\
\hline
& {\small (iv)} & {\small no} & \\
\hline
\end{tabular}

\medskip

\

\begin{center}
Table 2 
\end{center}

\smallskip

 \begin{tabular}{|c|c|c|c|}
\hline
\quad $\kk$ \quad & {\small Type} $t$ & {\small Existence results } & {\small \qquad Parameters of $t$  for the existence \qquad } \\ \hline
$\hh_1$ & {\small (i)} & {\small yes} & {\small $\nu\neq 0, \mu=-\frac{\eta^{2}}{\nu}$} \\ 
& & & {\small $\eta=\nu=0$} \\ \hline
 & {\small (ii)} & {\small no} &   \\ \hline
$\rr_3$ & {\small (i)} & {\small no } & \\ \hline
$\rr_3$ & {\small (ii)} & {\small no } & \\ 
 \hline  
	\end{tabular}

\begin{tabular}{|c|c|c|c|}
\hline
\quad $\kk$ \quad & {\small Type} $t$ & {\small Existence results } & {\small \qquad Parameters of $t$  for the existence \qquad } \\ \hline
$\rr_{3,\lambda}$ & {\small (i)} & {\small yes } & {\small $\lambda=0, \nu\neq 0, \varepsilon=0$}\\
& & & {\small $\lambda=0, \eta=1, \nu=0, \varepsilon =0 $} \\
& & & {\small $\lambda=0, \eta=-1, \nu=0, \varepsilon=0 $}  \\
& & & {\small $\lambda=-1, \eta=1, \mu=0 $}\\
& & & {\small $\lambda=-1, \eta=\frac{j_{43}\mu + j_{53}}{j_{53}^{2}}, \nu=-\frac{(j_{43}\mu+2j_{53})j_{43}}{j_{53}^{2}}$}\\
& & & {\small $\lambda=-1, \eta=-1, \nu=2\frac{j_{42}}{j_{52}}, \mu=0 $} \\ \hline
& {\small (ii)} & {\small yes } &  {\small $ \lambda=0, \nu\neq 0, \mu=-\frac{\eta(\eta-1)}{\nu}, \varepsilon=0$}\\
& & & {\small $ \lambda=0, \eta=1, \nu=0, \varepsilon=0$} \\
& & & {\small $ \lambda=0, \eta=\nu=0, \varepsilon=0$} \\
& & & {\small $ \lambda=0, \eta=\mu=0, \varepsilon=0$} \\
& & & {\small $ \lambda=0, \eta=\nu=0, \mu=-\frac{j_{53}}{j_{43}}, \varepsilon=0$} \\
& & & {\small $ \lambda=1, \nu\neq 0, \mu=-\frac{\eta(\eta-1)}{\nu}$} \\
& & & {\small $ \lambda=1, \eta=1, \nu=0$}\\
& & & {\small $ \lambda=1, \eta=\mu=0$} \\
& & & {\small $\lambda=1, \eta=\nu=0$} \\
& & & {\small $\lambda\neq 0, \eta=\frac{\mu j_{43}-\lambda j_{53}+j_{53}}{j_{53}}, \nu=-\frac{j_{43}(\mu j_{43}-2\lambda j_{53}+j_{53})}{j_{53}^{2}} $} \\
& & & {\small $\lambda\neq 0, \eta=\lambda, \mu=0$}\\
\hline
$\rr_{3,\delta}$ & {\small (i)} & {\small no} & {\small $$ }  \\ 
  \hline
  & {\small (ii)} & {\small yes} & {\small $\delta=0, \nu\neq 0, \mu=-\frac{\eta(\eta-1)}{\nu}$} \\ 
& & & {\small $\delta=0, \eta=1, \nu=0$} \\
& & & {\small $\delta=0, \eta=\mu=0 $} \\
& & & {\small $\delta=0, \eta=\nu=0$}\\
\hline
\end{tabular}

\medskip

The explicit matrix realizations of $j:\kk\to \RR^3$ and $j:\kk\to \hh_1$ can be seen in \cite{CCO}.

\


Acknoledgements. RCS acknowledges partial support from the research project MTM2010-18556 of the MICINN (Spain).

I. Cardoso and G. Ovando were partially supported by SCyT-UNR and CONICET.

\end{document}